\documentclass[11pt,a4paper,fleqn]{article}

\usepackage{amsmath, amssymb, amsthm}
\usepackage[top=2.5cm, bottom=2.5cm, left=1.7cm, right=1.7cm]{geometry}
\usepackage[colorlinks=true,urlcolor=blue,pdfstartview=FitH, linkcolor=blue,citecolor=blue,bookmarksopen, bookmarksnumbered={true}]{hyperref}
\usepackage{enumerate}
\usepackage{cancel}
\usepackage{mathpazo}

\usepackage{pdfsync}
\usepackage{color}
\newtheorem{theorem}{Theorem}[section]
\newtheorem{definition}{Definition}[section]
\newtheorem{corollary}{Corollary}[section]

\newtheorem{remark}{Remark}[section]

\definecolor{pink}{rgb}{1,0.08,0.58}
\definecolor{orange}{rgb}{1,0.5,0}
\definecolor{purple}{rgb}{0.75002,0,1}
\definecolor{olive}{RGB}{85,107,47}
\definecolor{mygreen}{rgb}{0,0.6,0}

\def\hh{{\mathcal{H}}}
\def\RR{{\mathbb{R}}}

\def\W{{\mathcal{W}}}
\def\B{{\mathcal{B}}}
\def\U{{\mathcal{U}}}
\def\V{{\mathcal{V}}}
\def\D{{\mathcal{D}}}
\date{}

\title{{\color{red}  Some consequences of Caristi's fixed point theorem, partial answers to some known open problems and its applications}}

\author{Farshid Khojasteh\footnote{Department of Mathematics, Arak-Branch, Islamic Azad University, Arak, Iran,\;\href{mailto:fr\_khojasteh@yahoo.com}{fr\_khojasteh@yahoo.com}.
                         }
\; , \; Erdal Karapinar\footnote{Department of Mathematics, Atilim University, $\dot{I}$ncek, 06836 Ankara, Turkey,\;\href{mailto:erdalkarapinar@yahoo.com}{erdalkarapinar@yahoo.com}.
}
\; , \; Hassan Khandani\footnote{Department of Mathematics, Mahabad-Branch, Azad University, Mahabad, Iran,\;\href{mailto:khandani.hassan@yahoo.com}{khandani.hassan@yahoo.com}.
}
}

\begin{document}
\maketitle
\date{}
\hrule

\abstract{
{\color{blue}
In this paper, we show that several extension of Banach contraction principle,
can be easily derived from the Caristi's theorem is one of the useful generalization of Banach contraction principle in the setting of the complete metric spaces. Moreover, some partial answers to some known open problems are given via Caristi's corollaries. Finally, existence of bounded solutions of a functional equation
is studied to support our results.
}
}
\newline
\\
{\bf\color{blue} Key words:} Caristi's fixed point, Hausdorff metric, Functional equation, Boyd and Wong's contraction, Meir-Keeler type.\\
{\bf\color{blue} 2010 MSC :} {\normalsize } 47H10,
54E05.\\
\hrule

{\color{red}
\section{Introduction and preliminaries}}
In the literature, the Caristi fixed-point theorem  is known as  one of the very interesting and useful generalization of the Banach fixed point theorem for self-mappings on a complete metric space. In fact,  Caristi fixed-point theorem is a modification of the $\varepsilon$-variational principle of Ekeland
(\cite{E1974, E1979}) that is a crucial tool in the nonlinear analysis, in particular, optimization, variational inequalities, differential equations.
and control theory. Furthermore, in 1977 Western  \cite{W1977}  proved that the conclusion of Caristi's theorem is equivalent to metric completeness.
In the last decades, Caristi's fixed-point theorem has been generalized and extended in several directions (see e.g., \cite{1,2}
and the related references therein).

{\par The Caristi's fixed point theorem asserted as follows:
\begin{theorem}\cite[Caristi]{CAR}\label{t1}
Let $(X,d)$ be a complete metric space and let $T:X\to X$ be a mapping such that
\begin{equation}\label{eq1}
d(x,Tx)\leq \varphi(x)-\varphi(Tx)
\end{equation}
where $\varphi:X\to[0,+\infty)$ be a lower semi continuous mapping. Then $T$ has at least a fixed point.
\end{theorem} }
Let us recall some basic notations, definitions and well-known results
needed in this paper. Throughout this paper, we denote by $\mathbb{N}$ and $\mathbb{R}$, the sets of positive integers and real numbers, respectively. Let $(X,d)$
be a metric space. Denote by $\mathcal{CB}(X)$ the family of all nonempty closed and
bounded subsets of $X$. A function $\mathcal{H}:\mathcal{CB}(X)\times
\mathcal{CB}(X)\rightarrow \lbrack 0,\infty )$ defined by
\begin{equation*}
\mathcal{H}(A,B)=\text{max}\left\{ \sup_{x\in B}d(x,A)\text{,}\sup_{x\in
A}d(x,B)\right\}
\end{equation*}
is said to be the Hausdorff metric on $\mathcal{CB}(X)$ induced by the
metric $d$ on $X$. A point $v$ in $X$ is a fixed point of a map $T$ if $v=Tv$
(when $T:X\rightarrow X$ is a single-valued map) or $v\in Tv$ (when $%
T:X\rightarrow \mathcal{CB}(X)$ is a multi-valued map).

In this work, we show that many of known Banach contraction's generalization can be deduce and generalize by Caristi's fixed point theorem and its consequences. Also, a partial answer to a known open problem is given via Caristi's corollaries.
{\color{red}\section{Main Result}}
In this section, we show that many of known fixed point results can be derived from Caristi's theorem.
{\color{blue}\begin{corollary}\label{co1}
Let $(X,d)$ be a complete metric space, and let $T:X\to X$ be a mapping such that
\begin{equation}
d(x,y)\leq \varphi(x,y)-\varphi(Tx,Ty),
\end{equation}
for all $x,y\in X$, where $\varphi:X\times X\to [0,\infty)$ is a lower semi continuous with respect to first variable. Then $T$ has a unique fixed point.
\end{corollary}}
\begin{proof}
For each $x\in X$, let $y=Tx$ and $\psi(x)=\varphi(x,Tx)$. Then for each $x\in X$
\[
d(x,Tx)\leq \psi(x)-\psi(Tx)
\]
and $\psi$ is a lower semi continuous mapping. Thus, applying Theorem \ref{t1} conclude desired result. To see the uniqueness of fixed point suppose that
$u,v$ be two distinct fixed point for $T$. Then
\[
d(u,v)\leq \varphi(u,v)-\varphi(Tu,Tv)=\varphi(u,v)-\varphi(u,v)=0
\]
Thus, $u=v$.
\end{proof}
{\color{blue}\begin{corollary}\cite[Banach contraction principle]{BAN}
Let $(X,d)$ be a complete metric space and let $T:X\to X$ be a mapping such that for some $\alpha\in[0,1)$
\begin{equation}\label{e6}
d(Tx,Ty)\leq \alpha d(x,y)
\end{equation}
for all $x,y\in X$. Then $T$ has a unique fixed point.
\end{corollary}}
\begin{proof}
Define, $\varphi(x,y)=\frac{d(x,y)}{1-\alpha}$. Then (\ref{e6}) shows that
\begin{equation}
(1-\alpha)d(x,y)\leq d(x,y)-d(Tx,Ty).
\end{equation}
It means that
\begin{equation}
d(x,y)\leq \frac{d(x,y)}{1-\alpha}-\frac{d(Tx,Ty)}{1-\alpha}
\end{equation}
and so
\begin{equation}
d(x,y)\leq \varphi(x,y)-\varphi(Tx,Ty)
\end{equation}
and so by applying Corollary \ref{co1}, one can conclude that $T$ has a unique fixed point.
\end{proof}
{\color{blue}\begin{corollary}
Let $(X,d)$ be a complete metric space and let $T:X\to X$ be a mapping such that
\begin{equation}\label{e8}
d(Tx,Ty)\leq \eta(d(x,y))
\end{equation}
where $\eta:[0,+\infty)\to[0,\infty)$ be a lower semi continuous mapping such that $\eta(t)<t$, for each $t>0$ and $\frac{\eta(t)}{t}$ be a non-decreasing map. Then $T$ has a unique fixed point.
\end{corollary}}
\begin{proof}
Define, $\varphi(x,y)=\frac{d(x,y)}{1-\frac{\eta(d(x,y))}{d(x,y)}}$, if $x\neq y$ and otherwise $\varphi(x,x)=0$. Then (\ref{e8}) shows that
\begin{equation}
(1-\frac{\eta(d(x,y))}{d(x,y)})d(x,y)\leq d(x,y)-d(Tx,Ty).
\end{equation}
It means that
\begin{equation}
d(x,y)\leq \frac{d(x,y)}{1-\frac{\eta(d(x,y))}{d(x,y)}}-\frac{d(Tx,Ty)}{1-\frac{\eta(d(x,y))}{d(x,y)}}
\end{equation}
Since $\frac{\eta(t)}{t}$ is non-decreasing and $d(Tx,Ty) < d(x,y)$ thus
\begin{equation}
d(x,y)\leq \frac{d(x,y)}{1-\frac{\eta(d(x,y))}{d(x,y)}}-\frac{d(Tx,Ty)}{1-\frac{\eta(d(Tx,Ty))}{d(Tx,Ty)}}=\varphi(x,y)-\varphi(Tx,Ty)
\end{equation}
and so by applying Corollary \ref{co1}, one can conclude that $T$ has a unique fixed point.
\end{proof}
Let $(X,d)$ be a complete metric space and $T :X\to X$ a map. Suppose there
exists a function $\phi: [0,+\infty) \to [0,+\infty)$ satisfying $\phi(0) = 0$, $\phi(s) < s$ for $s > 0$ and that
$\phi$ is right upper semi-continuous such that
\[
d(Tx,Ty)\leq\phi(d(x,y)) \ \ \ x,y\in X.
\]
Boyd-Wong \cite{Boyd} showed that $T$ has a unique fixed point.

Later, Meir-Keeler \cite{MKC} extended Boyd-Wong's result to mappings satisfying the
following more general condition:

For all $\epsilon > 0$ there exists $N > 0$ such that
\begin{equation}\label{mk}
 \epsilon\leq d(x,y) <\epsilon + \delta \ \ \ implies \ \ \ d(Tx,Ty)<\epsilon
\end{equation}
In 2001, Lim \cite{Lim} characterized condition (\ref{mk}) via introducing a category of functions namely L-functions which we recall his main results here.
\begin{definition}\cite{Suz,Lim}
Function $\phi$ from $[0,+\infty)$ into itself is called an $L-$function
if $\phi(0) = 0$, $\phi(s)>0$ for all $s\in(0,+\infty)$, and for every $s\in (0,+\infty)$ there exists $\delta>0$ such that $\phi(t)\leq s$, for all $t\in[\delta,\delta+\epsilon]$.
\end{definition}
In the case that $T:X\to \mathcal{CB}(X)$ is a multi-valued function Meir-Keeler condition can be written as:

For all $\epsilon > 0$ there exists $N > 0$ such that
\begin{equation}\label{mk}
 \epsilon\leq d(x,y) <\epsilon + \delta \ \ \ implies \ \ \ \hh(Tx,Ty)<\epsilon
\end{equation}
A known problem which is remains open until now, introduced by Lim in \cite{Lim} which is asserted as the following:\\
{\color{blue} Open Problem 1.}
\\
Let $(X,d)$ be a complete metric space and $T:X\to \mathcal{C}(X)$ a multi-valued
mapping such that $Tx$ is closed for every $x$ and
\[
\hh(Tx,Ty)<\phi(d(x,y)),
\]
for all $x,y\in X$, where $\phi$ is an $L-$function. Does $T$
have a fixed point?\\
The answer is yes if $Tx$ is compact for every $x$ (Reich \cite{reich}).

Another analogous open problem, raised in 2010 by Amini-Harandi\cite{Amini}.
In what follows, $\gamma:[0,+\infty)\to [0,+\infty)$ be subadditive, i.e.$\gamma(x+y)\leq \gamma(x)+\gamma(y)$, for each $x,y\in[0,+\infty)$, a nondecreasing continuous
map such that $\gamma^{-1}(\{0\})=\{0\}$, and let $\Gamma$ consist of all such functions. Also, let $\mathcal{A}$ be the class of all maps $\theta:[0,+\infty)\to [0,+\infty)$ for which there exists an $\epsilon_0 > 0$ such that
\[
\theta(t)\leq \epsilon_0 \ \ \Rightarrow \ \ \theta(t)\geq\gamma(t)
\]
where $\gamma\in\Gamma$.\\
{\color{blue} Open Problem 2.}\\
Assume that $T:X\to \mathcal{CB}(X)$ is a weakly contractive set-valued map on a complete metric space $(X,d)$, i.e.,
\[
\hh(Tx,Ty)\leq d(x,y)-\theta(d(x,y))
\]
for all $x,y\in X$, where $\theta\in\mathcal{A}$. Does $T$ have a fixed point?\\
The answer is yes if $Tx$ is compact for every $x$ (Amini-Harandi \cite[Theorem 3.3]{Amini}).

In the following theorem, we give partial answers to the above problems via Theorem \ref{co1} by replacing
some simple conditions instead of compactness condition $Tx$.\\
{\color{blue}\begin{theorem}\label{tt1} Let $(X,d)$ be a complete metric space, and let $T:X\to \mathcal{CB}(X)$ be a multi-valued function such that
\[
\hh(Tx,Ty)\leq \eta (d(x,y))
\]
for all $x,y \in X$, where $\eta :[0,\infty)\rightarrow [0,\infty)$ is lower semi continuous map such that $\eta(t)<t$, for all $t\in (0,+\infty)$ and $\frac {\eta(t)}{t}$ is non-decreasing.
Then $T$ has a fixed point.
\end{theorem}}
\begin{proof}Let $x\in X$ and and $y\in Tx$. If $y=x$  then $T$ has a fixed point and the proof is complete, so we suppose that $y\neq x$. Define
\[\theta (t)=\frac{\eta(t)+t}{2}\hbox{   for all $t\in (0,+\infty)$}\]
Since  $\hh(Tx,Ty)\leq \eta (d(x,y))<\theta (d(x,y))<d(x,y)$, so there exists $z\in Ty$ such that
\[d(y,z)<\theta (d(x,y))<d(x,y)\]
We again suppose that $y=z$, therefore $d(x,y)-\theta (d(x,y))\leq d(x,y)-d(y,z)$ or equivalently
\[d(x,y)<\frac{d(x,y)}{1-\frac{\theta (d(x,y))}{d(x,y)}}-\frac{d(y,z)}{1-\frac{\theta (d(x,y))}{d(x,y)}}\]
since $\frac{\theta(t)}{t}$ is also a nondecreasing function and $d(y,z)<d(x,y)$  we get
\[d(x,y)<\frac{d(x,y)}{1-\frac{\theta (d(x,y))}{d(x,y)}}-\frac{d(y,z)}{1-\frac{\theta (d(y,z))}{d(y,z)}}.\]
Define $\Phi(x,y)=\frac{d(x,y)}{1-\frac{\theta (d(x,y))}{d(x,y)}}$ if $x\neq y$, otherwise $0$ for all $x,y\in X$. Now arguing by induction we get a sequence $\{x_{n}\}$  such that $x_{n+1}\in Tx_{n}$ for all nonnegative integer $n$ and
\[d(x_{n},x_{n+1})\leq \Phi(x_{n},x_{n+1})- \Phi(x_{n+1},x_{n+2}).\]
Since $\{\Phi(x_{n},x_{n+1})\}$ is a decreasing and bounded from below sequence, so it converges. Also we deduce that $\{x_{n}\}$ is a cauchy sequence, because for all $m,n\in\mathbb{N}$ with $m>n$ we have
\[
\begin{array}{lll}
d(x_n,x_m)&\leq&\overset{m-1}{\underset{i=n}{\sum}}d(x_i,x_{i+1})\\\\
&\leq &\overset{m-1}{\underset{i=n}{\sum}}\Phi(x_{i},x_{i+1})- \Phi(x_{i+1},x_{i+2})\\\\
&=&\Phi(x_{n},x_{n+1})- \Phi(x_{m},x_{m+1})
\end{array}
\]
Thus,
\[
\lim_{n\to\infty}(\sup\{d(x_m,x_n):m>n\})=0.
\]
Therefore, by the completeness of $X$ it converges to some point $x\in X$. Now we show that $x$ is a fixed point of $T$.
\begin{equation}\label{e22}
\begin{array}{lll}
d(x,Tx)&\leq &d(x,x_{n+1})+d(x_{n+1},Tx)\\\\
&= &d(x,x_{n+1})+d(Tx_{n},Tx)\\\\
&\leq & d(x,x_{n+1})+\eta(d(x_{n},x)).
\end{array}
\end{equation}
By taking limit on both side of (\ref{e22}), we get $d(x,Tx)=0$ and this means that $x\in Tx$.
\end{proof}
{\color{blue}\begin{theorem}\label{tt2} Let $(X,d)$ be a complete metric space, and let $T:X\to \mathcal{CB}(X)$ be a multi-valued function such that
\[
\hh(Tx,Ty)\leq d(x,y)-\theta(d(x,y))
\]
for all $x,y \in X$, where $\theta :(0,\infty)\rightarrow (0,\infty)$ is lower semi continuous map such that, for all $t\in (0,+\infty)$ and $\frac {\theta(t)}{t}$ is non-increasing.
Then $T$ has a fixed point.
\end{theorem}}
\begin{proof}
Let $\eta(t)=t-\theta(t)$, for each $t>0$. Then, $\eta(t)<t$, for each $t>0$ and $\frac{\eta(t)}{t}=1-\frac{\theta(t)}{t}$ is non-decreasing. Thus, desired result is obtained by Theorem \ref{tt1}.
\end{proof}


{\color{blue}\begin{corollary}\cite[Mizoguchi-Takahashi's type]{Miz}
Let $(X,d)$ be a complete metric space and let $T:X\to X$ be a mapping such that
\begin{equation}\label{e7}
d(Tx,Ty)\leq \eta(d(x,y))d(x,y)
\end{equation}
where $\eta:[0,+\infty)\to[0,1)$ be a non-decreasing mapping. Then $T$ has a unique fixed point.
\end{corollary}}
\begin{proof}
Let $\theta (t)=\eta(t)t$ . $\theta(t)<t$ for all $t\in R_{+}$  and $\frac{\theta(t)}{t}=\eta (t)$  is a nondecreasing mapping. By the assumption $ d(Tx,Ty)\leq \eta (d(x,y))d(x,y)=\theta(d(x,y))$ for all $x,y\in X$, therefore by corollary 1.1 $T$ has unique fixed point.
\end{proof}
{\color{blue}\begin{corollary}\cite[Rhoades' type]{Rho}
Let $(X,d)$ be a complete metric space and let $T:X\to X$ be a mapping such that
\begin{equation}\label{e9}
d(Tx,Ty)\leq d(x,y)-\eta(d(x,y))
\end{equation}
where $\eta:[0,+\infty)\to(0,+\infty)$ be a mapping such that $\frac{\eta(t)}{t}$ be non-decreasing map. Then $T$ has a unique fixed point.
\end{corollary}}
\begin{proof}
Define $\theta(t)=t-\eta(t)$, by our assumption $\theta(t)<t$ for each $t> 0$ and $\frac{\theta(t)}{t}=1-\frac{\eta(t)}{t}$ is a non-decreasing map and $d(Tx,Ty)\leq d(x,y)-\eta(d(x,y))=\theta (d(x,y))$ , therefor the result follows by Theorem \ref{tt1}.
\end{proof}
{\color{red}\section{Existence of bounded solutions of functional equations}
}
Mathematical optimization is one of the fields in which the methods of fixed point theory are widely
used. It is well known that the dynamic programming provides useful tools for mathematical optimization
and computer programming. In this setting, the problem of dynamic programming related to multistage
process reduces to solving the functional equation
\begin{equation}\label{sw9}
p(x)=\sup_{y\in \D}\{f(x,y)+\Im(x,y,p(\eta(x,y)))\}, \ \ \ x\in \W,
\end{equation}
where $\eta:\W\times \D\to \W$, $f:\W\times \D\to \RR$ and $\Im:\W\times \D\times \RR\to \RR$.
We assume that $\U$ and $\V$ are Banach spaces,
$\W\subset \U$ is a state space and $\D\subset \V$ is a decision space. The studied process consists
of {\color{blue} a state space}, which is the set of the initial state, actions and transition model of the process and {\color{blue} a decision space}, which is the set of possible actions that are allowed for the process.

Here, we study the existence of the bounded solution of the functional equation \ref{sw9}.
Let $\B(\W)$ denote the set of all bounded real-valued functions on $W$ and, for an arbitrary $h\in \B(\W)$, define $||h|| = \sup_{x\in \W} |h(x)|$. Clearly, $(\B(W),||.||)$ endowed with the metric d defined by
\begin{equation}\label{p2}
d(h,k)=\sup_{x\in\W}|h(x)-k(x)|,
\end{equation}
for all $h,k\in \B(\W)$, is a Banach space. Indeed, the convergence in the space $\B(\W)$ with respect to $||.||$ is uniform. Thus, if we consider a Cauchy sequence $\{h_n\}$ in $\B(\W)$, then $\{h_n\}$ converges uniformly to a function, say $h^*$, that is bounded and so $h\in \B(\W)$.\\
We also define $T:\B(\W)\to \B(\W)$ by
\begin{equation}\label{re2}
T(h)(x)=\sup_{y\in \D}\{f(x,y)+\Im(x,y,h(\eta(x,y)))\}
\end{equation}
for all $h\in \B(\W)$ and $x\in \W$.

We will prove the following theorem.
{\color{blue}\begin{theorem}
Let $T : \B(\W)\to \B(\W)$ be an upper semi-continuous operator defined by (\ref{re2}) and assume that the
following conditions are satisfied:
\begin{itemize}
\item[$(i)$] $f:\W\times \D\to \RR$ and $\Im:\W\times \D\times \RR\to \RR$ are continuous and bounded;
\item[$(ii)$] for all $h,k\in \B(\W)$, if
\begin{equation}\label{pp2}
\begin{array}{rll}
0<d(h,k)<1 &  implies&  |\Im(x,y,h(x))-\Im(x,y,k(x))|\leq \frac{1}{2}d^2(h,k)\\
d(h,k)\geq 1 & implies & |\Im(x,y,h(x))-\Im(x,y,k(x))|\leq \frac{1}{2}d(h,k)
\end{array}
\end{equation}
where $x\in \W$ and $y\in \D$.

Then the functional equation (\ref{sw9}) has a bounded solution.
\end{itemize}
\end{theorem}
}
\begin{proof}
Note that $(\B(\W),d)$ is a complete metric space, where $d$ is the metric given by (\ref{p2}). Let $\mu$ be an
arbitrary positive number, $x\in W$ and $h_1, h_2 \in \B(\W)$, then there exist $y_1, y_2 \in \D$ such that
\begin{eqnarray}
T(h_1)(x)&<&f(x,y_1)+\Im(x,y_1,h_1(\eta(x,y_1)))+\mu,\\
T(h_2)(x)&<&f(x,y_2)+\Im(x,y_2,h_2(\eta(x,y_2)))+\mu,\\
T(h_1)(x)&\geq& f(x,y_1)+\Im(x,y_1,h_1(\eta(x,y_1))),\\
T(h_2)(x)&\geq&f(x,y_2)+\Im(x,y_2,h_2(\eta(x,y_2))).
\end{eqnarray}
Let $\varrho :[0,\infty)\rightarrow [0,\infty)$ be defined by
\[
\varrho(t)=\left\{
\begin{array}{ll}
\frac {1}{2}t^2, & \hbox{$0<t<1$} \\
\frac{1}{2}t, & \hbox{$t\ge 1$}
\end{array}
\right.
\]
Then we can say that (\ref{pp2}) is equivalent to
\begin{equation}\label{oo3}
|\Im(x,y,h(x))-\Im(x,y,k(x))|\leq\varrho(d(h,k))
\end{equation}
for all $h,k\in \B(\W)$. It is easy to see that $\eta(t)<t$, for all $t>0$ and $\frac{\eta(t)}{t}$ is a non-decreasing function.

Therefore, by using (19) and (22), it follows that
\[
\begin{array}{lll}
T(h_1)(x)-T(h_2)(x)&<& \Im(x,y_1,h_1(\eta(x,y_1)))-\Im(x,y_2,h_2(\eta(x,y_2)))+\mu\\
&\leq&|\Im(x,y_1,h_1(\eta(x,y_1)))-\Im(x,y_2,h_2(\eta(x,y_2)))|+\mu\\
&\leq&\varrho(d(h_1,h_2))+\mu.
\end{array}
\]
Then we get
\begin{equation}\label{gh3}
T(h_1)(x)-T(h_2)(x)<\varrho(d(h_1,h_2))+\mu.
\end{equation}
Analogously, by using (20) and (21), we have
\begin{equation}\label{gh4}
T(h_2)(x)-T(h_1)(x)<\varrho(d(h_1,h_2))+\mu.
\end{equation}
Hence, from (\ref{gh3}) and (\ref{gh4}) we obtain
\[
|T(h_2)(x)-T(h_1)(x)|<\varrho(d(h_1,h_2))+\mu,
\]
that is,
\[
d(T(h_1),T(h_2))<\varrho(d(h_1,h_2))+\mu.
\]
Since the above inequality does not depend on $x\in\W$ and $\mu > 0$ is taken arbitrary, then we conclude
immediately that
\[
d(T(h_1),T(h_2))\leq\varrho(d(h_1,h_2)),
\]
so we deduce that the operator $T$ is an $\varrho$-contraction. Thus, due to the continuity of $T$, Theorem \ref{tt1} applies to the operator $T$, which has a fixed point $h^*\in \B(\W)$, that is, $h^*$ is a bounded solution of the functional equation
(\ref{sw9}).
\end{proof}
\section*{Competing interests}
The authors declare that they have no competing interests.
\section*{Authors contributions}
All authors contributed equally and significantly in writing this
article. All authors read and approved the final manuscript.

\end{document}